\documentclass[reqno]{amsart}
\usepackage{amsmath, amssymb}

\newbox\Text
\newcount\width
\newcount\height
\newdimen\dimh
\newdimen\dimv
\newcount\dep
\newdimen\shifth
\newdimen\shiftv
\newcount\ellmarg

\newcount\cx
\newcount\cy

\newbox\O
\newdimen\dimhO
\newdimen\dimvO
\newcount\hO
\newcount\vO

\newbox\FText
\newdimen\marg
\newdimen\rs
\newdimen\scrspc

\def\encirc#1{
    \ellmarg=30         % size of the margin space inside ellipses
    \rs=1pt             % raise of the inscribed text over the baseline
    \scrspc=-0.5pt      % the value of \scriptspace
    \setbox\O=\hbox{$\bigcirc$}
    \dimhO=\wd\O
    \dimvO=\ht\O
    \hO=\dimhO
    \vO=\dimvO
    \divide\hO by 9472
    \divide\vO by 9472
    \advance\hO by \ellmarg
    \advance\vO by \ellmarg
    \advance\hO by -19
    \advance\vO by -10
    \setbox\Text=\hbox{$\scriptstyle{#1}\scriptspace=\scrspc$}
    \dimh=\wd\Text
    \dimv=\ht\Text
    \dep=\dp\Text
    \divide\dep by 9472
    \width=\dimh
    \divide\width by 9472
    \cx=-\width
    \advance\width by \ellmarg
    \height=\dimv
    \divide\height by 9472
    \cy=-\height
    \advance\cy by \dep
    \advance\height by \dep
    \advance\height by \ellmarg
    \advance\hO by -\width
    \advance\vO by -\height
    \ifnum\hO>0
        \ifnum\vO>0
            \putcircle
        \else
            \putellipse{\width}{\height}{\cx}{\cy}
        \fi
    \else
        \putellipse{\width}{\height}{\cx}{\cy}
    \fi
}

\def\putellipse#1#2#3#4{
\mskip8mu \raise\rs\hbox{%
\copy\Text%
\special{pn 5}%
\special{ar \the#3 \the#4 \the#1 \the#2 0 6.28}%
}\mskip8mu
}

\def\putcircle{
\shifth=\dimhO \advance\shifth by \dimh \divide\shifth by 2
\mskip3mu\hbox{\copy\O\kern-\the\shifth\raise\rs\hbox{\copy\Text}}\mskip7mu
}

\def\1{%
\hbox{1\kern-3.5pt 1\kern-4pt%
\raise 5.7pt\hbox{\vrule height1pt width1.5pt depth0pt}} }

\theoremstyle{plain}
\newtheorem{theorem}{Theorem}
\newtheorem{lemma}[theorem]{Lemma}

\newtheorem{corollary}[theorem]{Corollary}
\newtheorem{conjecture}[theorem]{Conjecture}

\theoremstyle{definition}
\newtheorem{definition}[theorem]{Definition}
\newtheorem{example}[theorem]{Example}

\theoremstyle{remark}
\newtheorem{remark}[theorem]{Remark}

\numberwithin{theorem}{section}
\numberwithin{equation}{section}
\newcommand{\thref}[1]{Theorem~{\rm\ref{#1}}}
\newcommand{\leref}[1]{Lemma~{\rm\ref{#1}}}
\newcommand{\coref}[1]{Corollary~{\rm\ref{#1}}}

\newcommand{\exref}[1]{Example~{\rm\ref{#1}}}

\newcommand\lbb[1]{\label{#1}}
\newcommand{\CC}{\Bbbk}
\newcommand{\ZZ}{\mathbb{Z}}
\renewcommand\geq{\geqslant}
\renewcommand\leq{\leqslant}
\DeclareMathOperator{\rk}{rk} \DeclareMathOperator{\Span}{Span}
\DeclareMathOperator{\End}{End} \DeclareMathOperator{\gk}{GKdim}
\DeclareMathOperator{\trdeg}{tr.deg} \DeclareMathOperator{\ad}{ad}
\DeclareMathOperator{\cur}{Cur} \DeclareMathOperator{\dif}{Dif{}f}
\DeclareMathOperator{\coeff}{Coef{}f}
\DeclareMathOperator{\cend}{Cend}
\DeclareMathOperator{\gc}{gc}\DeclareMathOperator{\Res}{Res}
\def\ti#1{\widetilde{#1}}
\def\wti{\widetilde}
\def\bar#1{\overline{#1}}
\def\d{\partial}
\def\c#1#2#3{#1\encirc{#2}#3}
\begin{document}

\title{On associative conformal algebras of linear growth II}
\author{Alexander Retakh}

\address{Department of Mathematics\\
Massachusetts Institute of Technology\\Cambridge, MA 02139}
\email{retakh@math.mit.edu}

\thanks{Paritally supported by the NSF}

\begin{abstract} We classify unital associative conformal algebras
of linear growth and provide new examples of such.
\end{abstract}

\maketitle

\section*{Introduction}\lbb{sec.intro}
Conformal algebras were introduced in \cite{K1} to provide
algebraic formalism for the singular part of the OPE in the theory of 
vertex algebras.  Since then
they turned out to be a useful instrument in the study of vertex
algebras (see, e.g., \cite{Ro}), infinite-dimensional Lie
superalgebras \cite{K3}, and in a generalized form, hamiltonian structures 
in integrable systems \cite{BDK}.

\begin{definition}\lbb{maindef}  A conformal algebra $C$ is a
$\CC[\d]$-module endowed with bilinear operations
$\encirc{n}:C\otimes C\to C$, $n\in\ZZ_{\geq 0}$ such that for any
$a,b\in C$
\begin{enumerate}
\item (locality axiom) $\c anb=0$ for $n>N(a,b)$

($N(a,b)$ is called the locality degree of $a$ and $b$);

\item $\d(\c anb)=\c{(\d a)}nb+\c an{(\d b)}$;

\item $\c{(\d a)}nb=-n\c a{n-1}b$.
\end{enumerate}
\end{definition}

The number $n$ in $\encirc{n}$ is called the order of
multiplication $\encirc{n}$.

In this paper $\CC$ is an algebraically closed field of zero
characteristic.  When clear, we refer to objects that are finite
as modules over $\CC[\d]$ simply as finite.

One of the basic questions in the study of conformal algebras is
the theory of representations of finite modules. Besides being of
independent interest, it is also related to the study of
representations of algebras of differential operators
\cite{BKL2,Ze2}.   Thus ones of the most important associative
conformal algebras are the algebras $\cend_n$ which are the
analogues of matrix algebras in ``ordinary'' theory.  From an
algebraic point of view, these are exactly simple unital
associative conformal algebras of linear growth \cite{Re1}.  The
next logical step is to describe all unital associative conformal
algebras of linear growth.  This paper contains such description.

\begin{theorem}\lbb{theo}  Let $C$ be a unital conformal algebra
of linear growth.  Then
\begin{itemize}
\item if $C$ is prime, then $C$ is isomorphic to
either $\cend_n$ or a subalgebra of a current algebra over a prime
algebra of linear growth;
\item if $C$ is semisimple, then $C$ embeds into a direct sum of
$\cend_n$ and a current algebra over a semiprime algebra of zero
or linear growth;
\item if $C$ contains a nilpotent ideal, then its coefficient
algebra is not semiprime.
\end{itemize}
\end{theorem}

The paper is virtually self-contained.  The first section contains
all necessary definitions and statements and basic examples of
associative conformal algebras.

The second section is devoted to classification results and contains the 
proof of \thref{theo} (see Theorems~\ref{confprime}, \ref{confsemisim}, 
and \leref{nilpid}). We deal
with prime conformal algebras first. The proofs follow along the
lines of those in \cite{Re1}; however, since we work in a more
general setup, we need to repeat some of the steps for consistency
of presentation.

The third section is concerned with subalgebras of prime conformal
algebras; in particular we exhibit a non-current unital subalgebra
of a current algebra.

{\bf Acknowledgements:} the results of this paper first appeared
in my thesis.  I am grateful to my advisor Efim Zelmanov.  The
text was typeset with the use of {\tt conformal.sty} package by
Michael Roitman.

\section{Preliminaries}\lbb{sec.prelim} A detailed exposition of
conformal algebras can be found in \cite[Chapter 2]{K1} and in the
survey papers \cite{K2,Ze1}.  The presentation here is guided by
our needs in the next sections and is in no way complete.  We also
present several results on unital conformal algebras from
\cite{Re1,Re2}.

Standard algebraic terminology easily carries over to the
conformal case, thus an {\em ideal} $I$ of a conformal algebra $C$ is a
conformal subalgebra such that $\c CnI\subset I, \c InC\subset I$
for all $n$, a {\em nilpotent} ideal $I$ is such that the product of a
fixed number of copies of $I$ (with multiplications of any order)
is zero, a {\em simple} conformal algebra contains no ideals, a
{\em semisimple} one contains no nilpotent ideals, etc.

\subsection{Basic examples}\lbb{subsec.maindef}
Let $A$ be any algebra.  We call elements of the extension
$A[[z,z^{-1}]]$ {\em formal distributions} on $A$.  Let
$f(z)=\sum_{n\in\ZZ} f(n)z^{-n-1}$ and $g(z)=\sum_{n\in\ZZ}
g(n)z^{-n-1}$ be formal distributions on $A$.  Define the product
$\encirc{m}$, $m\in\ZZ_{\geq 0}$, of $f(z)$ and $g(z)$ as
\begin{equation*}
\c {f(z)}m{g(z)}=\Res_{w=0} f(w)g(z)(w-z)^n
\end{equation*}
(by $\Res_{w=0}h(w,z)$ we mean a formal distribution in $z$ that
is a coefficient at $w^{-1}$ in $h(w,z)$ viewed as a formal
distribution on the set $A[[z,z^{-1}]]$).  Two formal
distributions are called local if only a finite number of such
products are non-zero.

A set of mutually local formal distributions that is closed with
respect to products $\encirc{m}$ and the operator $\d/\d z$ is a
conformal algebra.  Clearly, as the application of $\d/\d z$
preserves locality (though changes it degree), the $\CC[\d/\d
z]$-span of a set of mutually local distributions closed with
respect to products $\encirc{m}$ forms a conformal algebra.
Moreover, we have

\begin{lemma}[Dong's lemma \cite{Li,K1}]\lbb{donglemma} Let $f,g$,
and $h$ be pairwise mutually local formal distributions over
either a Lie or associative algebra. Then for any $n\geq 0$, $\c
fng$ and $h$ are again pairwise mutually local.
\end{lemma}

Thus, mutually local formal distributions over a Lie or
associative algebra generate a conformal algebra.

It can be shown \cite{K1} (see also \cite{Bo} for a vertex algebra
version of the same construction) that every conformal algebra $C$
embeds into formal distributions on some algebra.  Moreover, there
exists a universal algebra $\coeff C$ (called the {\em coefficient
algebra} of $C$) such that for any $A$, $C\to A[[z,z^{-1}]]$, there
exists a unique homomorphism $\coeff C\to A$ such that the diagram
\begin{equation*}
\begin{array}{c}
\coeff C[[z,z^{-1}]]\xrightarrow{\hspace{10pt}} A[[z,z^{-1}]]\\
\hspace{10pt}\text{\Large $\nwarrow$} \hspace{20pt} \text{\Large
$\nearrow$}\\[-2pt]
\hspace{10pt}C
\end{array}
\end{equation*}
commutes.  In particular, $\coeff C$ ``distinguishes''
coefficients $f(n)$ of every element of $f\in C$ (or rather the
corresponding formal distribution $f(z)\in\coeff $).  The subalgebra of 
coeffiecients at $z^{-1}$ is denoted $(\coeff C)_0$.

A conformal algebra is called {\em associative} (respectively, {\em Lie})  
if the corresponding coefficient algebra is associative (respectively,
Lie).  For every identity satisfied by $\coeff C$, one can write a
corresponding conformal identity satisfied by $C$;  however, for this
exposition we do not require the explicit forms of conformal
associativity, Jacobi identity, etc.

\begin{example}\lbb{ex.curr} Let $B$ be any algebra.  For every
$b\in B$ consider the following formal distribution on
$B[t,t^{-1}]$:
\begin{equation*}
\ti b=\sum_n bt^nz^{-n-1}.
\end{equation*}
Clearly for any $b_1,b_2\in B$, the formal distributions
$\wti{b_1}$ and $\wti{b_2}$ are mutually local:
$\c{\wti{b_1}}m{\wti{b_2}}=\delta_{0,m}\wti{b_1b_2}$.  Thus by
Dong's Lemma~\ref{donglemma}, $\ti b$ generate a conformal
algebra.  It is called the {\em current algebra} over $B$
and is denoted $\cur B$.

Observe that the conformal algebra generated by formal
distributions $\sum bz^{-n-1}$ on $B$ is isomorphic to $\cur B$.
However, for this conformal algebra $B$ is not a coefficient
algebra $\coeff\cur B$ , whereas $B[t,t^{-1}]$ is.
\end{example}

\begin{example}\lbb{ex.confweyl} Denote by $W$ the Weyl algebra
$\CC\langle x,t\,|\,xt-tx=1\rangle$ and by $W_t$ its localization
at $t$. We define the conformal algebra $\cend_n$ as an algebra of
formal distributions on $W_t$ generated by distributions
$L^k_A=\sum Ax^kt^nz^{-n-1}, k\geq 0, A\in\End_n(\CC)$.

In particular, $\cend_1$ is generated by elements $L^k=\sum
x^kt^nz^{-n-1}$ for $k=0,1$. Their non-zero products are
\begin{equation}\lbb{eq.confweylgen}
\begin{split}
&\c {L^0}0{L^0}=L^0,\quad \c {L^0}0{L^1}=\c {L^1}0{L^0}=L,\\
&\c {L^0}1{L^1}=\c {L^1}1{L^0}=-L^0,\\
&\c {L^1}0{L^1}=L^2,\quad \c {L^1}1{L^1}=-{L^1}.
\end{split}
\end{equation}
It follows that $\cend_1$ (and, as a consequence, $\cend_n$ for
any $n$) is not finite over $\CC[\d]$.

Observe also that $L^0$ generates a subalgebra isomorphic to
$\cur\CC$ and that, in a broad sense, $L^0$ acts as (left) identity.
This will be used later.

In the theory of representations of conformal algebras $\cend_n$
plays the role $\End_n(\CC)$ in ordinary theory, i.e., $\cend_n$
is the conformal algebra of conformal linear maps on $\CC[\d]^n$.
\end{example}

\begin{example}\lbb{ex.diff} Let $B$ be an associative algebra
with a locally nilpotent derivative $\delta$.  For every $b\in B$
consider the following formal distribution in
$B[t,t^{-1};\delta][[z,z^{-1}]]$:
\begin{equation*}
\ti b=\sum_n bt^nz^{-n-1}
\end{equation*}
(here $B[t,t^{-1};\delta]$ is the Ore extension of $B$ localized
at $t$, i.e., for any $b\in B$, $bt-tb=\delta(b)$).  Such formal
distributions are mutually local with each other, namely
\begin{equation}\lbb{eq.diff}
\c {\wti{b_1}}m{\wti{b_2}}=(-1)^m\wti{b_1\delta^m(b_2)},\quad
b_1,b_2\in B
\end{equation}
and locality follows from nilpotence of $\delta$.

The family of formal distributions $\{\ti b\,|\,b\in B\}$ spans
the conformal algebra called the {\em differential} conformal algebra
$\dif B$. Obviosuly, $(\coeff\dif B)_0=B$.

Observe that for a trivial $\delta$, $\dif B=\cur B$.  Also,
$\cend_n=\dif \CC[x]$ with the standard derivation.

It was shown in \cite{Re2} that $B$ and $\dif B$ have equivalent
categories of representations.  The lattice of ideals of $B$ is
isomorphic to the lattice of $\delta$-stable ideals of $\dif B$.
Namely, to a $\delta$-stable ideal $I$ of $B$ there corresponds
the ideal $\ti I$ of $\dif B$ spanned by $\{\ti b\,|\,b\in I\}$.
Conversely, to an ideal $J$ of $\dif B$ there corresponds the
ideal $\bar J=\{b\,|\,\ti b\in J\}$ of $B$ and $\wti{\bar J}=J$,
$\bar{\ti I}=I$.

As a corollary we have the following
\begin{lemma}\lbb{xcopsimple}
$\dif A$ is simple if and only if $A$ is differentiably simple.
\end{lemma}
\end{example}

\subsection{Unital conformal algebras}\lbb{subsec.unital}
The study of ordinary associative algebras begins with the study
of unital ones, i.e. those containing $\CC$.  By analogy, in the
conformal case one should start by considering associative
conformal algebras that contain a subalgebra of rank $1$ acting
faithfully in some sense.

It can be shown that $\cur\CC$ is the unique associative conformal
algebra with non-zero multiplication that is free of rank $1$.
Moreover, it can be shown that every module $M$ over $\cur\CC$
splits as $M=M_0\oplus M_1$, where $\c{\cur\CC}n{M_0}=0$ for every
$n$ and every element of $M_1$ is fixed by the action of
$e\encirc{0}$ for any generator $e$ of $\cur\CC$ \cite{Re2}.  This 
motivates the following definition:

\begin{definition}\lbb{def.unit} An associative conformal algebra
$C$ is {\em unital} if $\cur\CC$ embeds into $C$ and for the action of
the image of this embedding, $C=C_1$.
\end{definition}

A generator of $\cur\CC\subset C$ is called a {\em conformal identity}
and is denoted $e$.  Observe that a conformal identity is not
unique (see \leref{innerderiv}).

Unlike in the ordinary case, it is not clear how one can ``adjoin identity
'' to a torsion-free conformal algebra (a unital conformal algebra is
automatically torsion-free).  However, a differential conformal algebra
$\dif B$ always embeds into a unital conformal algebra: for this one needs
only to adjoin identity to $B$.  Thus we will always assume below that a
differential conformal algebra is unital.  Unless stated otherwise, we
will also assume that $e=\ti 1$.

Modulo a technical condition, the converse of the above
observation is also true.  For a conformal algebra $C$ denote by
$L(C)$ the set $\{a\,|\,\c anb=0 \,\forall b\in C,n\in\ZZ_{\geq
0}\}$.

\begin{theorem}\lbb{unitisdiff}  Let $C$ be a unital conformal
algebra such that $L(C)=0$.  Then $C=\dif B$ for a unital
associative algebra $B$.
\end{theorem}

In particular, a semisimple unital conformal algebra is always
differential.

\subsection{Gelfand-Kirillov dimension}\lbb{subsec.gkdim}

he Gelfand-Kirillov dimension of a finitely generated algebra (of
any variety) $A$ is defined as
\begin{equation*}
\gk A=\limsup_{r\to\infty}\frac{\log\dim(V^1+V^2+\dots+V^r)}{\log
r},
\end{equation*}
where $V$ is a generating subspace of $A$ \cite{KL}.  This
definition easily carries over to the conformal case.

Let $C$ be a finitely generated conformal algebra (over any
variety).  Define $C_r$ to be the $\CC[\d]$-span of products of at
most $r$ generators with any positioning of brackets and
multiplications of any order.

Since the powers of $\d$ can be gathered at the beginning of
conformal monomials (with a probable change in the orders of
multiplications), it is clear that $\bigcup_r C_r=C$.  For a given
ordered collection of generators and a given positioning of
brackets, the number of non-zero monomials is finite because of
locality.  Therefore, $\rk C_r$ is finite.

\begin{definition}\lbb{gkdimdef} Let $C$ be a finitely generated
conformal algebra. Then
\begin{equation}\lbb{eq.gkdimdef}
\gk C=\limsup_{r\to\infty}\frac{\log\rk_{\CC[\d]}C_r} {\log r}.
\end{equation}
\end{definition}

Conformal Gelfand-Kirillov dimension has the same basic properties
as the ordinary one: it is invariant of the choice of the
generating set, $\gk$ of a subalgebra or a quotient algebra does
not exceed that of the algebra, etc.

For a conformal associative algebra $C$, $\gk \coeff C\leq\gk C+1$
\cite[Theorem 2.2]{Re1} (the inequality is sometimes strict, e.g.,
when $C$ is torsion).  One can show directly that for a
differential conformal algebra, $\gk \dif B=\gk B$.  In
particular, $\gk \cend_n=1$.  The main result of \cite{Re1} is the
following

\begin{theorem}\lbb{confmainthclassif} Let $C$ be a simple unital
associative conformal algebra of Gelfand-Kirillov dimension $1$.
Then $C$ is isomorphic to $\cend_n$ for some $n$.
\end{theorem}

In the next section we present the generalization of this theorem.

\section{Semisimple conformal algebras of linear growth}\lbb{sec.semisim}
In this section we generalize \thref{confmainthclassif} and
achieve the complete generalization of unital conformal algebras
of $\gk$ $1$.

As follows from \leref{xcopsimple}, $C$ is simple if and only if
$\coeff C$ is differentiably simple.  Following this
correspondence, we define a larger subclass of unital associative
conformal algebras: we call $C$ {\em prime} whenever $(\coeff
C)_0$ is prime (since being prime is, in some sense, equivalent to
being differentiably prime). Also, recall that $C$ is {\em
semisimple} if it does not contain non-zero nilpotent ideals.  The
latter condition is equivalent to having a semiprime coefficient
algebra (see \leref{nilpid}).

\subsection{Classification of associative algebras of linear growth}  The
following theorem was proven in \cite{SSW} (see also \cite{SW}):

\begin{theorem}\label{ssw}  Let $A$ be a finitely generated algebra of
linear growth.  Then

(i) The nilradical N(A) of $A$ is nilpotent.

(ii) If $A$ is semiprime (i.e. if $N(A)=0$), then it is a finite
module over its center $Z(A)$ which is also finitely generated.
\end{theorem}

Several facts from the original proof of this theorem will be used
below as well.

\subsection{Prime unital conformal algebras of linear growth}  We are
going to classify prime unital conformal algebras of linear
growth.  It is well known that a differentiably simple algebra is
necessarily prime (see, e.g. \cite{Po}), hence
\thref{confmainthclassif} also follows from such a classification.

Let $A$ be a finitely generated prime algebra of linear growth.
Then by \thref{ssw}, it is a finite module over its center $Z(A)$.
Moreover, it is easy to see that for any derivation $\delta$ of
$A$, $Z(A)$ is $\delta$-stable: for $a\in Z(A)$,
$0=\delta([a,b])=[\delta(a),b]+[a,\delta(b)]=[\delta(a),b]$ for
any $b\in A$.

Thus, we begin by considering the case of a prime commutative
finitely generated algebra.

\begin{lemma}\lbb{commutprimecoeff} Let $A$ be a finitely generated prime
commutative algebra with a locally nilpotent derivation $\delta$,
$\gk A=1$.  Then either $A\cong\CC[x], \delta=\d/\d x$ or
$\delta=0$.
\end{lemma}

\begin{proof} Since $A$ is prime, a non-zero algebraic element of $A$ must
be invertible, hence all its algebraic elements lie in $\CC$.

Consider two sets of transcendental elements of $A$:
\begin{equation*}
\begin{aligned}
S_1&=\{x\in A\,|\,\text{all non-zero $\delta^n(x)$ are transcendental}\},\\
S_2&=\{x\in A\,|\,\text{for some $n$, $\delta^n(x)\neq 0$ and is
algebraic}\}.
\end{aligned}
\end{equation*}
Clearly, both sets are $\delta$-stable. Assume both are non-empty.
Without loss of generality we can pick $x_1\in S_1$, $x_2\in S_2$
such that $\delta(x_1)=0$, $\delta(x_2)=1$. As $\trdeg A=1$, $x_1$
and $x_2$ are algebraically dependent. In any statement of
dependence of $x_2$ over $\CC[x_1]$ the degree in $x_2$ can be
lowered by application of $\delta$.  Therefore, one of the $S_i$'s
is empty.

Consider now the case $A=\CC+S_1$.  Assume that there exists an
element with a non-zero derivation.  Without loss of generality we
can consider $x$ and $y$ such that $\delta(x)=y,\delta(y)=0$.
Just as above, consider a statement of dependence of $x$ over
$\CC[y]$.  Application of $\delta$ lowers the degree in $x$
(though it increases the degree in $y$), a contradiction.
Therefore, $\delta$ kills all transcendental elements.

The remaining case is $A=\CC+S_2$.  Choose $x$ such that
$\delta(x)=1$. Let $y$ be an arbitrary element with
$\delta(y)\in\CC$.  Then $x-y(\delta(y))^{-1}\in\CC$ and
$y\in\CC[x]$.  For an arbitrary $y\in S_2$, by induction on the
minimal $n$ such that $\delta^n(y)\in\CC$, we also obtain
$y\in\CC[x]$.
\end{proof}

\begin{remark}\lbb{wright} The final part of the proof of the above lemma
can be also deduced from a result in \cite{Wr}.
\end{remark}

\begin{corollary}\lbb{comuttdiffsimplecoeff} Let $A$ be a finitely
generated differentiably simple commutative algebra of growth $1$
with a locally nilpotent derivation.  Then $A\cong\CC[x]$.
\end{corollary}

\begin{proof} Indeed, if $A\not\cong\CC[x]$, it must be simple.
Therefore, $A$ is a field of transcendental degree $1$ and can not
be finitely generated.
\end{proof}

\begin{lemma}\lbb{primecoeff} Let $A$ be a finitely generated prime
algebra  with a locally nilpotent derivation $\delta$, $\gk A=1$.
Then $A$ is  either isomorphic to $\End_n(\CC[x]), \delta=\d/\d
x$, or $A$ can be embedded into an algebra $B$ such that $\delta$
extends to an inner derivation of $B$ determined by a nilpotent
element.
\end{lemma}

\begin{proof} As mentioned above, $A$ is a finite module over its center
$Z(A)$ which is finitely generated and has linear growth.
Clearly, $Z(A)$ is prime and $\delta$-stable.

{\bf Case 1.} $Z(A)=\CC[x]$ and $\delta|_{Z(A)}=\d/\d x$.

Consider subalgebra $A_0={\ker}\delta$.  We begin by demonstrating
that $A_0$ generates $A$ as a module over $Z(A)$. More precisely,
every $a\in A$ is of the form $\sum_1^n \dfrac{x^i}{i!}a_i$,
$a_i\in A_0$, where $n$ is such that $\delta^n(a)=0,
\delta^{n-1}(a)\neq 0$. Indeed, with the inductive assumption
$\delta(a)=\sum_1^{n-1} \dfrac{x^i}{i!} b_i$, $b_i\in A_0$,
consider $a_0=a-\sum_1^{n-1}\dfrac{x^{i+1}}{(i+1)!}b_i$. Since
$\delta(a_0)=0$, we see that $a$ is also a polynomial in $x$ over
$A_0$. Moreover, this polynomial form is unique for any $a\in A$.
Indeed, if $\sum_1^n x^ia_i=0, a_i\in A$, applying a necessary
number of derivations shows that the coefficient at the highest
power is $0$. In particular, this implies that $A=\CC[x]\otimes
A_0$.

Fix a subset $\{a_i\}$ of $A_0$ that generates $A$ as a module
over $Z(A)$. Any product of elements of $A_0$ is a linear
combination $\sum p_i(x)a_i$ over $Z(A)$ with the derivation $\sum
(\d p_i(x)/\d x)a_i=0$. This implies $A_0=\Span_{\CC}(a_i)$ is
finite dimensional.

Clearly any ideal of $A_0$ can be lifted to $A$, thus $A_0$ is
prime as well and therefore simple over $\CC$ \cite[2.1.15]{Rw}.
Hence, $A=\End_n(\CC[x])$ and $\delta=\d/\d x$.

{\bf Case 2.} $\delta|_{Z(A)}=0$.

Let $F$ be the field of fractions of $Z(A)$ and consider the
finite-dimensional simple $F$-algebra $B=F\otimes_{Z(A)}A$.
Clearly, $\delta$ extends to a derivation of $B$ and is,
therefore, an inner nilpotent derivation, $\delta=\ad a$
\cite{Ja}.  We may take $a$ to be nilpotent (it is enough to pick
any $a$ such that $\delta=\ad a$ and take its nilpotent part,
since the semisimple part must commute with all elements of $B$).
\end{proof}

\begin{remark}[\cite{Re1}]\lbb{proofmainth} If we strengthen the
condition of \leref{primecoeff} to $A$ being differentiably
simple, by \coref{comuttdiffsimplecoeff} we will have to consider
Case 1 of the above lemma only.  This implies
\thref{confmainthclassif}.
\end{remark}

\begin{remark}\lbb{conffingen} It follows that there exist no finitely
generated simple associative current conformal algebras of linear
growth. However, this (quite unexpected) result can be deduced
directly from \thref{ssw}.  Indeed, if $\cur A$ is such an
algebra, then $A$ must be simple and have linear growth.  So
should its center, hence it is a field of transcendental degree
$1$.
\end{remark}

In the more general framework of prime conformal algebras, the
second case of \leref{primecoeff} merits further consideration. We
notice first that by a change of conformal identity, one can
discount the twisting on the coefficient algebra introduced by an
inner derivation:

\begin{lemma}\lbb{innerderiv} Let $C$ be a differential conformal algebra
over algebra $A$ with an inner derivation determined by a
nilpotent element. Then $C\cong\cur A$. \end{lemma}

\begin{proof} We have $\delta=\ad r$ for a nilpotent $r$. Clearly,
$A[t,t^{-1};\delta]$ is isomorphic to the polynomial algebra
$A[s,s^{-1}]$ via the mapping $t-r\to s$.  To prove that
corresponding differential algebras $C$ and $\cur A$ are
isomorphic as well, we first show that the formal distribution
$e'=\sum (t-r)^n z^{-n-1}$ belongs to $C$.  Indeed, in this case,
the conformal subalgebra of $C$ generated by $e'$ and $\c {\ti
a}0{e'}, a\in A$, is isomorphic to $\cur A$.

Let $m$ be the degree of nilpotency of $r$. Since $\delta(r)=0$,
$t$ and $r$ commute; therefore,
\begin{equation*}
e'=\sum_{k=0}^{m-1} \frac{1}{k!} \d^k \left(\sum
r^kt^nz^{-n-1}\right)\,\in\,\dif A.
\end{equation*}

Conversely, $\c {(\c {\ti r}0{e'})}0{e'}=\ti r$ lies in the above
subalgebra, hence, so does $e$. Thus, this subalgebra coincides
with $C$.
\end{proof}

\begin{corollary}\lbb{confsimplefin} Let $C$ be a simple unital
associative conformal algebra that is finite over $\CC[\d]$.  Then
$C\cong\cur\End_n(\CC)$.
\end{corollary}

\begin{proof} Let $C=\dif A$ where $A$ is differentiably simple.  Since
$A$ is finite, it must be simple \cite{Bl}; thus
$A\cong\End_n(\CC)$ and all its derivations are inner.  The rest
follows from \leref{innerderiv}.
\end{proof}

\begin{remark}\lbb{confsimplefinrem} The above result also follows from
the classification of simple Lie conformal algebras that are
finite over $\CC[\d]$ in \cite{DK}.  In fact, in this line of
proof one does not require unitality; though, we still get only
$\cur\End_n(\CC)$ as an answer \cite{K2}.  This shows that every
simple associative conformal algebra that is finite over $\CC[\d]$
is unital.

Of course, in general not every associative conformal algebra is
unital and one can not simply ``adjoin'' a conformal identity as
in the ordinary case.  However, in every known case, a conformal
algebra can be embedded into a unital one.

\begin{conjecture}  Every $\d$-torsion free associative conformal algebra
can be embedded into a unital conformal algebra.
\end{conjecture}

Another question is: if such embeddings $C\to C'$ exist for a
given conformal algebra $C$, what is the lower bound on $\gk C'$?

Clearly, when $C$ has a faithful representation that is finite
over $\CC[\d]$, $C$ embeds into $\cend_n$.  Moreover, when $C$ is
finite itself it can be embedded into a finite conformal algebra
(this follows from the classification and so far no direct proof
is known).  So, in such cases it is always possible to find $C'$
with $\gk C'=\gk C$.
\end{remark}

We can now translate the statement of \leref{primecoeff} into the
language of conformal algebras:

\begin{theorem}\lbb{confprime} Let $C$ be a prime unital finitely
generated associative conformal algebra with $\gk C=1$.  Then $C$
is isomorphic to either $\cend_n$ or a subalgebra of a current
algebra over a finitely generated prime algebra.
\end{theorem}

\begin{proof} By \leref{primecoeff}, $(\coeff C)_0$ is either
$\End_n(\CC)\otimes\CC[x]$ with a nilpotent derivation given by
$\d/\d x$ or a subalgebra of a prime algebra with an inner
derivation determined by a nilpotent element. In the first case,
$C$ is isomorphic to $\cend_n$ and in the second case it is a
subalgebra of a current algebra by \leref{innerderiv}.  In
general, this current algebra might be infinitely generated but by
construction one may choose an appropriate prime current
subalgebra.
\end{proof}

Essentially, this classification is the best one can hope for;
this is explained in Section~\ref{sec.confprimeex}.

\subsection{Classification of semisimple conformal algebras of linear
growth}

We begin by translating semisimplicity of conformal algebras of
linear growth into a condition for its coefficient algebra.

\begin{lemma}\lbb{nilpid}  Let $\dif A$ be an associative conformal
algebra of zero or linear growth.  Then $\dif A$ is semisimple if
and only if $A$ is semiprime.
\end{lemma}

\begin{proof} By construction from \leref{xcopsimple}, an
ideal J of $\dif A$ is nilpotent if and only if the corresponding
ideal $\bar J$ of $A$ is nilpotent.  Indeed, if for any
$a_0,\dots,a_n$,
$(\dots(\wti{a_0}\encirc{0}\dots)\encirc{0}\wti{a_n}=0$, then
$a_0\cdot\ldots\cdot a_n=0$ and, conversely, if $(\bar J)^n=0$,
then $(\dots(\wti{a_0}\encirc{m_1}\dots)\encirc{m_n}\wti{a_n}=0$
as $\bar J$ is $\delta$--stable.

Thus, if $\dif A$ contains a nilpotent ideal, so does $A$.

Conversely, if $A$ is not semiprime, its nilradical $N(A)$ is
$\delta$-stable \cite[2.6.28]{Rw} and nilpotent (\thref{ssw}, (i)
for the case of linear growth).  Hence, $\wti{N(A)}$ is nilpotent.
\end{proof}

Thus, we are able to classify semisimple conformal algebras of
linear and zero growth.

We need an easy (and probably known) technical lemma first.

\begin{lemma}\lbb{derivrestrict} Let $A=\bigoplus_i A_i$ be a finite
direct sum of unital associative algebras $A_i$.  Then a
derivation $\delta$ on $A$ restricts to each $A_i$.
\end{lemma}

\begin{proof} Let $e_i$ be the identity in $A_i$.  Since
$\delta(e_ie_j)= e_i\delta(e_j)+\delta(e_i)e_j=0$ and the summands
lie in $A_i$ and $A_j$ respectively, we have $\delta(e_i)e_j=0$
(for any $j\neq i$).  Thus, $\delta(e_i)\in A_i$.  Now, as
$\delta(e_i)= \delta(e_i^2)= 2\delta(e_i)$, it follows that
$\delta(e_i)=0$. Consequently, $\delta(A_i)= \delta(A_ie_i)=
\delta(A_i)e_i\subset A_i$.
\end{proof}

The classification of semisimple unital conformal algebras of
$\gk\leq 1$ follows:

\begin{lemma}[\cite{K2}]\lbb{confsemisimfin} Let $C$ be a semisimple
unital finitely generated associative conformal algebra that is
finite over $\CC[\d]$. Then $C\cong\bigoplus_{i=1}^k
\cur\End_{n_i}(\CC)$.
\end{lemma}

\begin{proof} The proof is the same as for \coref{confsimplefin}. Let
$C=\dif A$.  By \leref{nilpid}, we have $A\cong \bigoplus_{i=1}^k
\End_{n_i}(\CC)$.  The nilpotent derivation that leads to $C$ is
inner, hence its effects can be removed by a change of conformal
identity as in \leref{innerderiv}.
\end{proof}

\begin{theorem}\lbb{confsemisim} Let $C$ be a semisimple unital finitely
generated associative conformal algebra, $\gk C= 1$. Then $C$
embeds into a direct sum of a current algebra over a semiprime
algebra of zero or linear growth and $\cend_n$ for some $n$.
\end{theorem}

\begin{proof} Let $C=\dif A$ be determined by a nilpotent derivation
$\delta$.  By \leref{nilpid}, A is semiprime.  It follows from
various lemmas in \cite{SSW} that $A$ splits as $A=B\oplus F$
where $B$ is semiprime Goldie and $F$ finite-dimensional.

We have $Q=Q(B)=\bigoplus_i Q_i$, a semisimple Artinian quotient
algebra of $B$.  By \leref{derivrestrict}, $\delta$ restricts to
$B$.  The standard construction of $Q(B)$ implies that $\delta$
can be extended to $Q$ and we can again restrict it to $Q_i$
(though it is no longer locally nilpotent at this stage).  We
obtain a system of prime ideals $P_i=B\cap \bigoplus_{i\neq j}
Q_j$.  Clearly, $P_i$ is $\delta$-stable, hence so is $B/P_i$.  As
$B\hookrightarrow \bigoplus_i B/P_i$ with the action of $\delta$
preserved, we have $\dif B\hookrightarrow\bigoplus_i \dif B/P_i$.
Also, $\dif A=\dif B\oplus \dif F$.

Thus, $\dif A$ embeds into a direct sum of prime conformal
algebras of linear growth and a $\CC[\d]$-finite semisimple
conformal algebra. We have two types of components in this sum:
some come from subalgebras of $\cend_m$, hence their sum can be
viewed as a subalgebra of $\cend_n$ for some $n$.  Others are
subalgebras of prime current algebras of growth not exceeding $1$.
Thus, their sum is a subalgebra of a semiprime current algebra.
\end{proof}

\section{Examples of subalgebras of prime unital conformal
algebras}\lbb{sec.confprimeex} In this section we discuss what
kinds of conformal algebras can live inside typical examples of
prime conformal algebras.

The most innocently looking one is a current algebra over a prime
algebra.

In the proof of \leref{innerderiv} we relied on the simple fact
that to each locally nilpotent derivation of $\coeff C$
corresponds a conformal identity and a canonical basis of $C$,
hence the effect of the derivation could be ``untwisted.''  This
is not necessarily true for all subalgebras of $C$.

\begin{remark}\lbb{sameidentity}  Let $C'$ be a unital conformal
subalgebra of a unital current algebra $C$ with the same conformal
identity. Then $C'$ itself is current.  Indeed, let
$a=\sum_0^n\d^k \wti{a_k}\in C'$.  Then $(-1)^n n!\wti{a_n}=\c
ane\in C'$, hence, by induction all $\wti{a_k}\in C'$.
\end{remark}

However, when conformal identities of the conformal algebra and
its subalgebra are different, the approach in the above remark can
not be used.  Such a situation arises in the setting of
\thref{confprime}: let $C'$ be a subalgebra of $C=\dif A$, where
the derivation on $A$ is inner, $\delta=\ad a$.  If $\ti a\not\in
C'$, the change of conformal identity prescribed by
\leref{innerderiv} can not be performed inside $C'$, thus $C'$
becomes a possibly non-current subalgebra of a current algebra
$C$.  Consider the following

\begin{example}\lbb{ex.noncur} Let $A=\End_2(\CC[x])$ and $\delta=\ad
e_{12}$ be a locally nilpotent derivation on $A$.  Remark that
$\dif A$ is current by \leref{innerderiv}. Let $B=\End_2(x\CC[x])$
with the identity adjoined.  $B$ is $\delta$-stable, so $\dif
B\subset\dif A$; however, we will show below that $\dif B$ is not
current for any choice of conformal identity.
\end{example}

More generally, it turns out that whenever the derivation is
external with respect to $C'$, $C'$ can not be current.  The
following statement is, in some sense, the converse of
\leref{innerderiv}.

\begin{lemma}\lbb{noncurlemma} Let $C$ be an associative conformal algebra
such that for different choices of conformal identities,
$C\cong\cur A$ and $C\cong\dif B$ for a nilpotent derivation
$\delta$.  Then $\delta$ can be made inner on $A$, i.e., there
exists $a\in A$ such that the conformal algebra $\dif A$
determined by $\ad a$ is isomorphic to $\dif B$ with the
isomorphism preserving the conformal identity.  Moreover, $a$ is
nilpotent.
\end{lemma}

\begin{proof}  We fix the following notations: $\ti a$, $a\in A$, is the
canonical basis of $\cur A$ and $\ti e$ is its conformal identity.
For $\dif B$, $\bar b$ is the canonical basis and $\bar e$ the
conformal identity.

We also identify the elements of $\cur A$ and $\dif B$ via the
given isomorphism.  Thus, we have $\bar e=\sum\d^i \wti{e_i}$ for
some $e_i\in A$. For an arbitrary $b\in B$, $\bar b=\sum\d^i
\wti{b_i}$.  Since $\bar b=\c {\bar b}0{\bar e}$, we have $\bar
b=\c {(\sum \d^i \wti{b_i})}0{(\sum \d^i \wti{e_i})}=\sum\d^i
\wti{b_0e_i}$. Thus, $b_i=b_0e_i$.

In particular this implies that if $b\neq 0$, then $b_0\neq 0$.
Moreover, $(bb')_0=b_0b'_0$ as $\overline{bb'}=\c {\bar
b}0{\overline{b'}}$.  This establishes a map $B\to A$, $b\mapsto
b_0$.  Now we have to show that it extends to a map of given
conformal algebras.

We also remark that as $\bar b=\c {\bar e}0{\bar b}$, we see that
$e_0b_0=b_0$.

According to (\ref{eq.diff}), $\overline{\delta(b)}=-\c {\bar
e}1{\bar b}$ for any $b\in B$.  We will now calculate $\c {\bar
e}1{\bar b}$ in $\cur A$:
\begin{equation}\lbb{eq.bar1}
\begin{split}
\c {\bar e}1{\bar b}&=
\c {(\sum \d^i\wti{e_i})}1{(\sum\d^i\wti{b_0e_i})}\\
&=\c {\wti{e_0}}1{(\sum\d^i\wti{b_0e_i})}-
\c {\wti{e_1}}0{(\sum\d^i\wti{b_0e_i})}\\
&=\sum \d^i(\c {\wti{e_0}}1{\wti{b_0e_i}})+ \sum i\d^{i-1}(\c
{\wti{e_0}}0{\wti{b_0e_i}})-
\sum\d^i\c {\wti{e_1}}0{\wti{b_0e_i}}\\
&=\sum i\d^{i-1}\wti{e_0b_0e_i}-\sum\d^i\wti{e_1b_0e_i}
\end{split}
\end{equation}
(the last equality is valid, as we are working with a current
algebra: products of positive orders of basis elements are $0$).

On the other hand,
\begin{equation}\lbb{eq.bar2}
\overline{\delta(b)}=\sum\d^i\wti{\delta(b)_0e_i}.
\end{equation}

Comparing the coefficients at $\d^0$ in (\ref{eq.bar1}) and
(\ref{eq.bar2}), we see that
$\delta(b)_0e_0=-e_0b_0e_1+e_1b_0e_0$.  Since
$\delta(b)_0e_0=\delta(b)_0$, we see that
$\delta(b)_0=\ad(e_1)b_0$.  In particular, this implies that
$\ad(e_1)$ is nilpotent.

Consider now a differential algebra $\dif A$ over $A$ determined
by $\ad(e_1)$.  By the above, we obtained a injective map
$(B,\delta)\to (A,\ad(e_1))$ of  differential algebras; therefore,
$\dif B$ embeds into $\dif A$ with the embedding given by $\bar
b\mapsto \ti{b_0}$.

It remains to show that such an embedding is surjective.  Since
$\cur A\cong\dif A$, it is possible to express $\ti a$, $a\in A$,
as $\ti a=\sum_j \d^j \overline{a_j}$, where $a_j\in B$.  Hence,
$\ti a=\sum_{i,j} \d^{i+j} \wti{(\overline{a_j})_i}$.  As $C$ is
free over $\CC[\d]$, $a=(\overline{a_0})_0$. This show that the
constructed map is an isomorphism.

For our purposes, we need also to show that $e_1$ is nilpotent.
Substitute $\bar b=\bar e$ in (\ref{eq.bar1}).  The coefficient at
$\d^j$ in the last line is $\wti{e_{j+1}}-\wti{e_1e_j}$. Since $\c
{\bar e}1{\bar e}=0$, we obtain by induction that $e_j=(e_1)^j$
for $j\geq 1$.  As the expression for $\bar e$, $\sum \d^i
\wti{e_i}$, is a finite sum, we see that $e_1$ is nilpotent.  This
completes the proof.
\end{proof}

\begin{corollary}\lbb{noncur} Let $A'$ be a subalgebra of $A$ and $a\in A$
a nilpotent element.  Then the conformal algebra $\dif A'$
determined by $\ad a$ is current if and only if there exists
$a'\in A'$, $\ad a=\ad a'$.
\end{corollary}

\begin{proof} If such $a'$ exists, the statement follows from
the proof of \leref{innerderiv}.

Otherwise, assume that $\dif A'$ is current for some choice of
conformal identity.  By \leref{noncurlemma}, $\dif A'$ is
determined by an inner derivation (for the same choice of
conformal identity). Hence, there exists such $a'$.
\end{proof}

This shows that the subalgebra in \exref{ex.noncur} is indeed
never current.

\begin{remark}\lbb{noncursubalg} Such examples can be constructed with
ease.  In particular, this means that the classification in
\thref{confprime} is the best possible for the case of prime
conformal algebras.
\end{remark}

Just as subalgebras of current algebras, subalgebras of $\cend_n$
also appear naturally.  However, this happens either when such a
subalgebra acts on a given finite module or in the context of
\thref{confsemisim}. In either case the conformal identities of
$\cend_n$ and its subalgebra coincide.  Thus, we essentially speak
of subalgebras of $\End_n(\CC)\otimes\CC[x]$ (with the same
identity).

Such subalgebras are too general to describe (cf. \thref{ssw}),
even in the prime case \cite{SW}.  The only known result is that
on differentiably simple subalgebras, i.e.
\thref{confmainthclassif} and \coref{confsimplefin}, and here the
resulting conformal subalgebras are isomorphic to either $\cend_m$
or $\cur\End_m(\CC)$.

The case of simple subalgebras is obviously important for
conformal representation theory.  Kac's conjecture \cite{K2}
describes all subalgebras of $\cend_n$ that act irreducibly on the
standard module $\CC[\d]^n$.  The unital conformal algebras from
the list are the ones described above; therefore, we can say that
our results in this section confirm Kac's conjecture for unital
algebras (see also \cite{BKL1}).

%%%%%%%%%%%%%%%%%%%%%%%%%%

%%%%%%%%%%%%%%%%%%%%%
\end{document}